\documentclass{article}
\usepackage[utf8]{inputenc}
\usepackage{amssymb,amsmath,latexsym}
\usepackage[english]{babel}

\usepackage{mathtools}

\title{Field Generated by Division Points of Certain Formal Group Laws - III}
\author{Soumyadip Sahu\footnote{Date : October25, 2019}}
\date{}
\begin{document}
\maketitle
\begin{abstract}
Several questions about the Galois group of field generated by certain one dimensional formal group laws are studied. This is a continuation of treatment in \cite{some results}, \cite{some more results}.\\
General formalism related to $\pi$-divisible groups is discussed in appendix.   
\end{abstract}
\section{Introduction}
In \cite{some results} and \cite{some more results} author studied some questions about the Galois group of field generated by division points of certain formal group laws and relation between this Galois group and the ring of endomorphisms of  corresponding formal group law. Present article is a continuation of \cite{some more results} .\\~\\
First the set-up and main results are reviewed : \\~\\
Let $p$ be a prime and let $K$ be a finite extension of $\mathbb{Q}_p$. Put $O_K$ to be the ring of integers of $K$, let $\mathfrak{p}_K$  
denote the unique maximal ideal of $O_K$ and $v_{\mathfrak{p}_K}(\cdot)$ be the valuation associated to it. Fix an algebraic closure $\overline{\mathbb{Q}}_p$ and $|.|_p$ be an fixed extension of the absolute value. Let $\overline{O}$ be the ring of integers of $\overline{\mathbb{Q}}_p$ and $\overline{\mathfrak{p}}$ be the unique maximal ideal of $\overline{O}$. Clearly $\overline{\mathfrak{p}} \cap K = \mathfrak{p}_K$. For $n \in \mathbb{N}$ let $\mu_{n}$ denote the group of $n$-th roots of unity inside $\overline{\mathbb{Q}}_p$.\footnote{In our convention $0 \notin \mathbb{N}$.}\\~\\
Let $A$ be an integrally closed, complete subring of $O_K$. Assume that $\mathfrak{p}_A$ is the maximal ideal, $K(A)$ is the field of fractions and $k(A)$ is the field of residues. Put $[k(A) : \mathbb{F}_p] = f_A$. Let $\mathfrak{F}$ be a (one dimensional, commutative) formal group-law defined over $O_K$ admitting an $A$ module structure. $\pi$ be a generator of $\mathfrak{p}_A$ and say \[[\pi](X) = \pi X + a_2X^2 + a_3X^3 + \cdots \in O_K[[X]] \tag{1.1}\] with at least one $a_i \in O_K - \mathfrak{p}_K$. Then $\min \,\{ i \,|\, |a_i|_p = 1 \} = p^h$ for some positive integer $h$ (see \cite[18.3.1]{haz}). Now if $\pi_1$ is another generator of $\mathfrak{p}_A$ and \[[\pi_1](X) = \pi_1 X + b_2X^2 + \cdots \in O_K[[X]]\] then $b_{p^h} \in O_K - \mathfrak{p}_K$ and $\min \,\{ i \,|\, |b_i|_p = 1\} = p^h$. This integer $h$ is called the height of $\mathfrak{F}$ as formal $A$ module. If $a_i \in \mathfrak{p}_K$ for all $i \geq 2$ then we say, height of $\mathfrak{F}$ is infinity. We shall only consider formal $A$ modules of finite height.\\~\\
Let $A$ be as above and $\mathfrak{F}$ be a formal $A$ module over $O_K$. It defines a $A$ module structure on $\mathfrak{p}_K$ which naturally extends to a $A$-structure on $\overline{\mathfrak{p}}$. We shall denote the corresponding addition by $\oplus_{\mathfrak{F}}$ to distinguish it from usual addition.\\
Use $\text{End}_{O_K}(\mathfrak{F})$ to denote ring of endomorphisms of formal group $\mathfrak{F}$ which are defined over $O_K$ and $\text{End}_{O_K}^{A}(\mathfrak{F})$ to denote ring of endomorphisms of formal $A$ module $\mathfrak{F}$ which are defined over $O_K$. It is well-known that \[\text{End}_{O_K}(\mathfrak{F}) = \text{End}_{O_K}^{A}(\mathfrak{F})\;\; (\text{see}\,\cite[21.1.4]{haz}).\]\\
Let $\pi$ be a generator of $\mathfrak{p}_A$. For each $n \geq 1$, use $\mathfrak{F}[\pi^n]$ to denote the $\pi^n$-torsion submodule of $\overline{\mathfrak{p}}$.\\ For any sub-field $L$ of $\overline{\mathbb{Q}}_p$, let $L_{\mathfrak{F}}(\pi^n)$ be the subfield of $\overline{\mathbb{Q}}_p$ generated by $\mathfrak{F}[\pi^n]$ over $L$ and put 
\[ \begin{split}
\Lambda_{\pi}(\mathfrak{F}) = \bigcup_{i \geq 1} \mathfrak{F}[\pi^i],\\
L_{\mathfrak{F}}(\pi^{\infty}) = \bigcup_{i \geq 1} L_{\mathfrak{F}}(\pi^i).
\end{split}\]We shall adopt the convention $\mathfrak{F}[\pi^0] = \{0\}$.\\~\\
Let $\pi$ be a generator  of $\mathfrak{p}_K$ and put $K_{\pi} = \mathbb{Q}_p(\pi)$. Use $A_{\pi}$ to denote the ring of integers of $K_{\pi}$. For simplicity we shall write $\mathfrak{p}_{A_\mathfrak{\pi}}$ as $\mathfrak{p}_{\pi}$ and $f_{A_\pi}$ as $f_{\pi}$. Note that $\pi$ is a generator of $\mathfrak{p}_\pi$. Let $\mathfrak{F}$ be a formal $A_{\pi}$ module of height $h$ defined over $O_K$.\\
$A$ be an integrally closed subring of $O_K$ containing $A_{\pi}$ such that $\mathfrak{F}$ has a $A$ module structure. Clearly, $\mathfrak{F}$ has height $h$ as $A$ module. We know that $h_A\,|\,f$ (See \cite[Remark-1.3(v)]{some results}). Put $h_{r,A} = \frac {h} {f_A}$. $h_{r, A_{\pi}}$ will be abbreviated as $h_{r, \pi}$ and if $A$ is clear from context we shall drop $A$ from subscript.\\~\\
\textbf{Notations and terminologies :}\\
Let $K, O_K, \mathfrak{p}_K, \pi$ be as before. We shall use the following notations and terminologies --\\~\\
A local field is a finite extension of $\mathbb{Q}_p$,\\
$O_k$ = ring of integers of $k$ for any discrete valuation field $k$ ,\\
$\mathfrak{p}_k$ = the unique maximal ideal of $k$,\\
$E$ = a finite unramified extension of $K$,\\
$K^{(n)}$ = the unique unramified extension of $K$ of degree $n$ for any $n \in \mathbb{N}$, $O_K^{(n)}$ is the ring of integers of $K^{(n)}$, \\ 
$L = E^{\text{ur}} = K^{\text{ur}}$,\\
$\widehat{L}$ = completion of $L$,\\
$\mathcal{F}_{\pi}(O_E)$ = set of all $A_\pi$ modules of finite height defined over $O_E$,\\
$\mathcal{F}_{\pi}(O_E)(h)$ = set of all $A_{\pi}$ modules of height $h$ defined over $O_E$ for any $h \in \mathbb{N}$,\\
$\mathcal{F}_{\pi}(O_{\widehat{L}})$ = set of all $A_{\pi}$ modules of finite height defined over $O_{\widehat{L}}$,\\
$\mathcal{F}_{\pi}(O_{\widehat{L}})(h)$ = set of all $A_{\pi}$ modules of height $h$ defined over $O_{\widehat{L}}$ for any $h \in \mathbb{N}$,\\
$C_p$ = completion of $\overline{\mathbb{Q}}_p$,\\
$\mathcal{E}_{E}$ = set of all subextensions of $\overline{\mathbb{Q}}_p/E$,\\
$\mathcal{E}_{L}$ = set of all subextensions of $\overline{\mathbb{Q}}_p/L$,\\
$\mathcal{E}_{\widehat{L}}$ = set of all subextensions of $C_p/\widehat{L}$,\\
$G_{E} = \text{Gal}(\overline{\mathbb{Q}}_p|E)$,\\
$G_{L} = \text{Gal}(\overline{\mathbb{Q}}_p|L)$,\\
$G_{\widehat{L}} = \text{Gal}(C_p|\widehat{L})$.\\~\\
Note that $G_{L} = G_{\widehat{L}}$ and it is a closed normal subgroup of $G_E$, known as \emph{inertia subgroup}.\\
Let $\pi$ be a generator of $\mathfrak{p}_K$ . An $A_{\pi}$ module $\mathfrak{F}$ over $O_K$ is called \emph{$\pi$-unramified group law} over $O_K$. If $K/\mathbb{Q}_p$ is unramified one can take $\pi =p$ and any formal group law is $p$-unramified.\\
Let $\mathfrak{F} \in \mathcal{F}_{\pi}(O_E)(h)$ and use $F$ to denote fraction field denote fraction field of $\text{End}_{O_E}(\mathfrak{F})$. Clearly, $K_{\pi} \subseteq F \subseteq E$. One can show that $K_{\pi} \subseteq F \subseteq K_{\pi}^{(h_{r, \pi})}$ (see \cite[remark-4.1.3]{some results}). If $F = K_{\pi}^{(h_{r, \pi})}$ then we say $\mathfrak{F}$ has \emph{endomorphism ring with full height}. Such $\mathfrak{F}$ defines a full group law in Lubin's sense (see \cite[4.3.1]{lubin1}).\\~\\
One has the maps 
\[\begin{split}
D_{\pi, E} : \mathcal{F}_{\pi}(O_E) \to \mathcal{E}_E , \, \mathfrak{F} \to E_{\mathfrak{F}}(\pi^{\infty}) \\
D_{\pi}^{\text{ur}} : \mathcal{F}_{\pi}(O_L) \to \mathcal{E}_L , \, \mathfrak{F} \to L_{\mathfrak{F}}(\pi^{\infty})\\
\widehat{D}_{\pi} : \mathcal{F}_{\pi}(O_{\widehat{L}}) \to \mathcal{E}_{\widehat{L}}, \, \mathfrak{F} \to \widehat{L}_{\mathfrak{F}}(\pi^{\infty}).
\end{split}\]
If $E$ is clear from context we shall drop it from subscript.\\
Additionally, one needs to consider formal $\mathbb{Z}_p$ modules over $O_E$ and associated module of $p$-torsion points. Here height shall mean height as $\mathbb{Z}_p$ module. While working with this set-up we shall use similar notations with $\pi$ replaced by $p$.\\~\\
In section-2 and section-3 notion of $p$-adic lie groups over local field is used. We shall follow definitions and terminologies of \cite{serre 1}. \\~\\
The following remark summarizes the main points of development so far :\\~\\
\textbf{Remark 1.1 :} i) Let $\mathfrak{F} \in \mathcal{F}_{p}(O_{\widehat{L}})$. Then $\text{End}_{O_{\widehat{L}}}(\mathfrak{F}) = \text{End}(\mathfrak{F})$ where $\text{End}(\mathfrak{F})$ is the absolute endomorphism ring as in \cite[2.3.3]{lubin1}.\\  
ii) Let $\mathfrak{F} \in \mathcal{F}_{\pi}(O_E)$. Then \\
a) $\text{End}_{O_E}(\mathfrak{F}$), $\text{End}(\mathfrak{F})$ are complete, integrally closed $A_{\pi}$ algebras of finite type.\\
b) $A_{\pi} \subseteq \text{End}(\mathfrak{F}) \subseteq A_{\pi}^{(h_{r, \pi})}$.\\
c) Further if $\mu_{p^h - 1} \subseteq O_E$, $\text{End}_{O_E}(\mathfrak{F}) = \text{End}(\mathfrak{F})$. \\See \cite[Section-2]{some more results}.\\
iii) Let $\mathfrak{F} \in \mathcal{F}_{\pi}(O_E)(h)$ and assume that $\mu_{p^h - 1} \subseteq O_E$. Then the following are equivalent :\\
a) $\text{Gal}(D_{\pi}(\mathfrak{F})|E)$ is abelian.\\
b) For all $n \in \mathbb{N}$ and $z \in \mathfrak{F}[\pi^n] - \mathfrak{F}[\pi^{n-1}]$ we have $E(\pi^n) = E(z)$.\\
c) The endomorphism ring of $\mathfrak{F}$ over $O_E$ has full height ie $\text{End}_{O_E}(\mathfrak{F}) = A_{\pi}^{(h_{r,\pi})}$. \\
Statement $(b)$ can be replaced by :\\
For all $n \in \mathbb{N}$ there exists $z \in \mathfrak{F}[\pi^n] - \mathfrak{F}[\pi^{n-1}]$ such that $E(\pi^n) = E(z)$.\\
See \cite[Theorem-2.3]{some results} and \cite[Remark-4.1]{some more results}.\\  
iv) Hypothesis as in (iii). If the equivalent conditions hold then :\\
a) $D_{\pi}(\mathfrak{F})/E$ is totally ramified and the Galois group is isomorphic to $U_{A_{\pi}^{(h_{r,\pi})}}$, the group of units of $A_{\pi}^{(h_{r,\pi})}$. See \cite[Appendix-A]{the}.\\
b) $D_{\pi}^{\text{ur}}(\mathfrak{F}) = E^{\text{ab}}$. See \cite[Section-2]{some more results}.\\~\\
In \cite{some results}, \cite{some more results} the Galois representation on associated Tate module was studied. In this article we shall continue these kind of investigations, fill up a few gaps and generalize the set-up to lie groups over general local fields and discuss question about number of generators as in \cite[Section-4]{some more results}. \\
There is an appendix devoted to discuss formalism of $\pi$-divisible groups suitable for this purpose.\\~\\
\textbf{Remark 1.2 :} i) Let $R$ be a complete, discrete valuation ring of characteristic $(0,p)$ with perfect residue field. Divisible formal group-laws over $R$ (see appendix, section-B) form the category of so-called `connected $p$-divisible groups' over $R$. Note that in this situation, \'etale points of a $p$-divisible groups over $R$ are defined over an unramified extension of fraction field of $R$. Since we shall be concerned about representation of inertia subgroup, these results can be stated in terms of any $p$-divisible group of dimension 1 ie without loss of generality one can ignore \'etale points. \\
ii) Author was ignorant of some important results of Serre and Serre-Sen while writing \cite{some results} and \cite{some more results}. These are put into proper place in the development of ideas (see section-2 and section-3). This answers some questions posed in these articles which constitutes section-4.  

\section{Galois representation - I}
Let $\mathfrak{F} \in \mathcal{F}_p(O_E)(h)$. It is known that $p$-adic Tate module $T_p(\mathfrak{F})$ is a free $\mathbb{Z}_p$ module of rank $h$. Put $V_p(\mathfrak{F}) = T_p(\mathfrak{F})\underset{\mathbb{Z}_p}{\otimes} \mathbb{Q}_p$. Fixing an ordered base for $T_{p}(\mathfrak{F})$ one obtains continuous representations    
\[\begin{split}
\rho_{p}(\mathfrak{F}) : G_E \to \text{Gl}_{h}(\mathbb{Q}_p),\\
\widehat{\rho}_{p}(\mathfrak{F}) : G_{\widehat{L}} \to \text{Gl}_h(\mathbb{Q}_p).
\end{split}\]\\
Note that, these maps factor through canonical inclusion $\text{Gl}_h(\mathbb{Z}_p) \to \text{Gl}_h(\mathbb{Q}_p)$. The representations $\rho_p(\mathfrak{F})$ and $\widehat{\rho}_p(\mathfrak{F})$ are quite similar, we shall mostly be concerned with $\widehat{\rho}_p(\mathfrak{F})$ and mention $\rho_p(\mathfrak{F})$ only if needed ie we want to study the totally ramified part of the extension. Note that here we are considering $\mathfrak{F}$ to be a group law over $O_{\widehat{L}}$.\\~\\
Let $H$ be image of $\widehat{\rho}_p(\mathfrak{F})$. Clearly $H$ is a closed, compact $p$-adic lie sub-group of $\text{Gl}_h(\mathbb{Q}_p)$ contained in $\text{Gl}_h(\mathbb{Z}_p)$. Further, $\rho_p(\mathfrak{F})$ induces a topological isomorphism $\text{Gal}(\widehat{D}_p(\mathfrak{F})|\widehat{L}) \cong H$. Let $\mathfrak{h}$ be the associated lie algebra.\\
 Further $H_{\text{alg}}$ be the smallest algebraic subgroup of $\text{Gl}_h(\mathbb{Q}_p)$ containing $H$ and $\mathfrak{h}_{\text{alg}} \subseteq \mathfrak{gl}(h, \mathbb{Q}_p)$ be the smallest algebraic lie sub-algebra of $\mathfrak{gl}(h, \mathbb{Q}_p)$ containing $\mathfrak{h}$. One can verify $\mathfrak{h}_{\text{alg}}$ is indeed the lie algebra associated to $H_{\text{alg}}$.  The following results are due to Serre :\\~\\
\textbf{Proposition 2.1 :} Assume that the following holds :\\
i) $V_p(\mathfrak{F})$ is a semi-simple $H$ module.\\
ii) $\text{End}(\mathfrak{F}) = \mathbb{Z}_p$.\\
Then, $H_{\text{alg}} = \text{Gl}_h$ and $H$ is an open sub-group of $\text{Gl}_h(\mathbb{Q}_p)$ in metric topology.\\~\\
\textbf{Proof :} See \cite[Section-5, Theorem-4]{serre 2} . $\square$\\~\\
Now one would like to know when condition-(i) of proposition-2.1 holds. In this direction we have the following result :\\~\\
\textbf{Proposition 2.2 :} The following are equivalent :\\
i) $V_{p}(\mathfrak{F})$ is a semi-simple $H$ module.\\
ii) $V_{p}(\mathfrak{F})$ is a semi-simple $\mathfrak{h}$ module.\\
iii) $\mathfrak{h}$ is a reductive Lie algebra ie a product of an abelian and a semi-simple Lie algebra and $V_p(\mathfrak{F})$ is a semisimple $\mathfrak{c}$ module where $\mathfrak{c}$ is center of $\mathfrak{h}$.\\~\\
\textbf{Proof :} See \cite[Section-1, Proposition-1]{serre 2}. $\square$\\~\\
\textbf{Proposition 2.3 :} $V_{p}(\mathfrak{F})$ is a simple $\mathfrak{h}$ module .\\~\\
\textbf{Proof :} See \cite[Section-5, Proposition-8]{serre 2}. $\square$ \\~\\
As a consequence, we have $H$ is open subgroup of $\text{Gl}_{h}(\mathbb{Q}_p)$ provided $\text{End}(\mathfrak{F}) = \mathbb{Z}_p$. Note that, this means $H$ has finite index in $\text{Gl}_{h}(\mathbb{Z}_p)$.\\~\\
Before main results of the section, we recall another result due to Serre (\cite[Section-3, Thm. 9]{serre 3}) :\\~\\
\textbf{Proposition 2.4 :} Let $Y$ be a closed, smooth, analytic sub-manifold of $(\mathbb{Z}_p)^N$  and let $Y_n$ denote its image under reduction modulo $p^n$. Let $d$ be dimension of $Y$ and $\text{vol}(Y)$ be volume of $Y$ in the induced measure. Then $|Y_n| = \text{vol}(Y)p^{nd}$ for sufficiently large $n$.\\~\\
Note that $\text{Gl}_{h}(\mathbb{Z}_p)$ is a closed subset of $(\mathbb{Z}_p)^{h^2}$. Hence $H$ is closed sub-manifold of $(\mathbb{Z}_p)^{h^2}$.  Let $d$ be dimension of $H$. Define $\widetilde{H} = H - \text{Id}$ where $\text{Id}$ is identity matrix. Clearly it is also a closed, analytic sub-manifold of $(\mathbb{Z}_p)^{h^2}$ of dimension $d$ and $|H_n| = |\widetilde{H}_n|$ for all $n \geq 1$.\\~\\
For $n \in \mathbb{N}$, let $\text{red}_n : H \to \text{Gl}_{h}(\mathbb{Z}/p^n\mathbb{Z})$ be the canonical reduction modulo $p^n$ map. Put $K_n = \text{Ker}(\text{red}_n)$ and $I_n = \text{Im}(\text{red}_n)$. We have, $K_n \cong \text{Gal}(\widehat{L}_{\mathfrak{F}}(p^\infty)|\widehat{L}_{\mathfrak{F}}(p^n))$ and $I_n \cong \text{Gal}(\widehat{L}_{\mathfrak{F}}(p^n)|\widehat{L})$. Further $|I_n| = |\widetilde{H}_n|= \text{vol}(H)p^{nd}$ for large enough $n$.\\~\\
This shows that dimension of image plays important role in determining the degree of the extensions which in term has influence of determining arithmetical properties of the corresponding extension. In particular, note that $|\text{Gal}(\widehat{L}_{\mathfrak{F}}(p^{n+1})|\widehat{L}_{\mathfrak{F}}(p^n))| = |I_{n+1}|/|I_{n}| = p^{d}$  for large enough $n$. This is the main theme for rest of this section.\\~\\
Put $\mathfrak{F}'[p^n] = \mathfrak{F}[p^n] - \mathfrak{F}[p^{n-1}]$ for all $n \in \mathbb{N}$. Let $m_p(n)$ be the smallest size of a subset of $\mathfrak{F}'[p^n]$ which generate $\widehat{L}_{\mathfrak{F}}(p^n)$ over $\widehat{L}$. Fix such a set $\{x_{i}(n) \,|\, 1 \leq i \leq m_p(n)\}$ for each $n$. Note that if $\{z_1, \cdots, z_h\}$ is a $\mathbb{Z}_p$ base for $T_p(\mathfrak{F})$ then $\{z_i(n) \,|\, 1 \leq h\}$ generates $\widehat{L}_{\mathfrak{F}}(p^n)$ over $\widehat{L}$. Hence $m_p(n) \leq h$. \\~\\
Let $k, n \in \mathbb{N}$ with $ 1 \leq k \leq n$. Put $G_{n,k} = \text{Gal}(\widehat{L}_{\mathfrak{F}}(p^n)|\widehat{L}_{\mathfrak{F}}(p^k))$. Clearly, $|G_{n,k}| = |I_n|/|I_k|$.\\~\\  
\textbf{Lemma 2.5 :} i) The set theoretic map $\Delta_{n,k} : G_{n,k} \to (\mathfrak{F}[p^{n-k}])^{m_p(n)}$ defined by $\Delta_{n,k}(\sigma) = (\sigma(x_i(n)) \ominus_{\mathfrak{F}} x_i(n) )$ is injective.\\
ii) $\Delta_{n,k}$ is homomorphism of groups if $n - k \leq k$.\\
iii) $d \leq hm_p(n)$ for large enough $n$.\\~\\
\textbf{Proof :} Proof of $(i), (ii)$ is clear. $(iii)$ follows from taking $k = n - 1$ and using the estimate mentioned before. $\square$\\~\\
\textbf{Lemma 2.6 :} i) Let $n_0 \in \mathbb{N}$ be such that $n_0 \leq k$ implies $d = hm_p(k)$ and $|\widetilde{H}_k| = \text{vol}(H) p^{kd}$. Then for all $n_0 \leq k \leq n$, $\Delta_{n,k}$ is onto.\\
ii) Further if $n_0 \leq k\leq n \leq 2k$ \[G_{n,k} \cong (\mathfrak{F}[p^{n-k}])^{m_p(n)}.\] In particular, the fields $\widehat{L}(x_i(n), \mathfrak{F}[p^k])$ are linearly disjoint over $\widehat{L}_{\mathfrak{F}}(p^k)$.
\\~\\
\textbf{Proof :} i) Follows from the fact that \[|G_{n,k}|= |I_n|/|I_k|= p^{d(n-k)}= p^{hm_p(n)(n-k)} = |\mathfrak{F}[p^{n-k}]|^{m_p(n)}.\]\\
ii) First part follows from part $(i)$ and part $(ii)$ of lemma 2.5.\\
Second part follows from first part. $\square$ \\~\\
\textbf{Remark 2.7 :} i) If $d = h^2$ then by lemma-2.5(iii) $m_p(n) = h$ for large enough $n$ and hypothesis of corollary-2.6 is verified. Note that, this is the case if $H$ is open in $\text{Gl}_{h}(\mathbb{Z}_p)$. \\
ii) In Sen's theory of ramification in $p$-adic lie filtration one obtains similar estimates but with respect to a lie filtration which may be different from the filtration $\{K_n\}$ (see \cite{sen}). But one can show by analyzing trivilization maps of a neighbourhood of identity in $H$, that $\{K_n\}$ is indeed equivalent to a lie filtration in sense of \cite[Chapter 6]{probst} and hence one can use the ramification estimates derived in Sen's theorem. \\~\\
Since dimension of the lie group $H$ is an important invariant in determining the Galois group of $\widehat{D}_p(\mathfrak{F})/\widehat{L}$ one would like to know when two formal group of same height has same dimension of image of $\widehat{\rho}_p$. We conclude the section with a preliminary observation in this direction. \\~\\
Let $\mathfrak{F}, \mathfrak{G} \in \mathcal{F}_p(O_{\widehat{L}})(h)$ be such that $\text{Hom}_{O_{\widehat{L}}}(\mathfrak{F}, \mathfrak{G}) \neq 0$. We know $\widehat{D}_p(\mathfrak{F}) = \widehat{D}_p(\mathfrak{G})$ (see \cite[Remark-2.2]{some more results}). Put $H(\mathfrak{F}) = \text{Im}(\widehat{\rho}_p(\mathfrak{F}))$ and $H(\mathfrak{G}) = \text{Im}(\widehat{\rho}_p(\mathfrak{G}))$.  Define $\{I_n(\mathfrak{F})\}$ and $\{I_n(\mathfrak{G})\}$ as before. \\~\\
\textbf{Lemma 2.8 :} Notation be as above. Let $\text{dim}(H(\mathfrak{F})) = d_1$ and $\text{dim}(H(\mathfrak{G})) = d_2$. Then $d_1 = d_2$. \\~\\
\textbf{Proof :} In light of discussions in this section it is enough to show that the sequence $\{\frac {|I_n(\mathfrak{F})|} {|I_n(\mathfrak{G})|} \}$ is bounded above and does not converge to zero.\\
Let $0 \neq f \in \text{Hom}_{O_{\widehat{L}}}(\mathfrak{F}, \mathfrak{G})$. Note that there exists $g \in \text{Hom}_{O_{\widehat{L}}}(\mathfrak{G}, \mathfrak{F})$ with $g \neq 0$ and the induced morphisms 
\[\begin{split}
f : (\overline{\mathfrak{p}}, \oplus_{\mathfrak{F}}) \to (\overline{\mathfrak{p}}, \oplus_{\mathfrak{G}}) \\
g : (\overline{\mathfrak{p}}, \oplus_{\mathfrak{G}}) \to (\overline{\mathfrak{p}}, \oplus_{\mathfrak{F}})
\end{split}\]
has finite kernels and are surjective (see \cite[1.6]{lubin2}).\\
Further they induce surjective morphisms 
\[\begin{split}
f : \Lambda_p(\mathfrak{F}) \to \Lambda_p(\mathfrak{G})\\
g : \Lambda_p(\mathfrak{G}) \to \Lambda_p(\mathfrak{G}).
\end{split}\]   
Let $c \in \mathbb{N}$ be such that $[p^c]$ annihilates $\text{ker}(f)$ and $\text{ker}(g)$. Then $\mathfrak{G}[p^n] \subseteq f(\mathfrak{F}[p^{n+c}])$ and $\mathfrak{F}[p^n] \subseteq g(\mathfrak{G}[p^{n+c}])$ for all $n \in \mathbb{N}$. Thus $\widehat{L}_{\mathfrak{F}}(p^n) \subseteq \widehat{L}_{\mathfrak{G}}(p^{n+c})$ and $\widehat{L}_{\mathfrak{G}}(p^n) \subseteq \widehat{L}_{\mathfrak{F}}(p^{n+c})$. Now using the estimates from proposition-2.4 one concludes the result. $\square$

\section{Galois Representation - II}
Let $\mathfrak{F} \in \mathcal{F}_{p}(O_{\widehat{L}})$ and $A$ be an integrally closed, complete ring with $A \subseteq \text{End}(\mathfrak{F}) = \text{End}_{O_{\widehat{L}}}(\mathfrak{F})$. Let $\omega \in \mathfrak{p}_A$ be a generator and $\mathfrak{F}$ be a formal $A$ module of height $h_A$. Note that this means $\mathfrak{F}$ has height $eh_A$ as $\mathbb{Z}_p$ module where $e$ is the ramification index of $K(A)/\mathbb{Q}_p$. Write $h_p$ for $\mathbb{Z}_p$ height. Clearly, $h_p = eh_A$.\\ $T_{\omega}(\mathfrak{F}) = \varprojlim \mathfrak{F}[\omega^n]$ is a free $A$ module of rank $h_{r,A} = \frac {h_A} {f}$ where $f$ is the degree of residue extension of $K(A)$ (see \cite[Lemma-1.1]{lubin2} and \cite[Section-2]{the}).\\
Consider the canonical representation \[ \widehat{\rho}_{\omega}(\mathfrak{F}) : G_{\widehat{L}} \to \text{Gl}_A(T_{\omega}(\mathfrak{F}))\]
Denote the image by $H_{\omega}$. Clearly we have a topological isomorphism $H_{\omega} \cong \text{Gal}(\widehat{L}_{\mathfrak{F}}(\omega^{\infty})|\widehat{L}) = \text{Gal}(\widehat{L}_{\mathfrak{F}}(p^{\infty})|\widehat{L}) \cong H$ where $H$ is as in previous section.\\
Note that if we fix an ordered $A$-base for $T_{\omega}(\mathfrak{F})$, $\text{Gl}_A(T_{\omega}(\mathfrak{F}))$ can be identified with $\text{Gl}(h_{r,A},A)$ and $H_{\omega}$ is a closed sub-group of $\text{Gl}(h_{r,A},A)$. A natural question that can be asked if $H_{\omega}$ is a analytic lie group over $K(A)$, fraction field of $A$.\\
This question can be rephrased as follows :\\~\\
$T_p(\mathfrak{F})$ is a free $\mathbb{Z}_p$ module of rank $h_p$. Fixing an ordered base for $T_p(\mathfrak{F})$ we have a canonical $\mathbb{Z}_p$ algebra homomorphism \[A \to M(h_p, \mathbb{Z}_p)\]
where $M(h_p, \mathbb{Z}_p)$ is the algebra of $h_p\times h_p$ matrices with entries in $\mathbb{Z}_p$. Further, note that this homomorphism is injective (see \cite[Section-2.2]{tate}). We shall abuse notation and identify $A$ with its image. Image of $a \in A$ shall be denoted by $[a]$. Clearly $A$ commutes with $H$.\\
Let $\mathfrak{h} \subseteq \mathfrak{gl}(h_p,\mathbb{Q}_p)$ be the corresponding lie algebra. One would like to show $\mathfrak{h}$ has $A$ module structure. To be precise we want to show that there exists an open subgroup  $U \subseteq \mathfrak{h}$ such that $AU \subseteq U$. Note if this is the case then a neighbourhood of identity in $H\;(\cong H_{\omega}$ as topological groups) has $A$-structure and it can be endowed with a structure of $K(A)$ manifold.\\
Consider the extension $K(A)/\mathbb{Q}_p$ and let $K_0$ be the corresponding maximal unramified extension. Let $A_0$ be its ring of integers. We have, $A = A_0[\omega]$. Clearly, $AU \subseteq U$ if and only if $A_0U \subseteq U$ and $[\omega]U \subseteq U$.\\~\\
\textbf{Remark 3.1 :} Without loss of generality one can assume  $U \subseteq M(h_p,\mathbb{Z}_p)$ and it is a $\mathbb{Z}_p$ submodule.\\~\\   
Let $\{z_{1}, \cdots, z_{h_{r,A}}\}$ be a $A$-base for $T_{\omega}(\mathfrak{F})$. Clearly $\{z_1(n), \cdots z_{h_{r,A}}(n)\}$ generates $\widehat{L}_{\mathfrak{F}}(\omega^n)$ over $\widehat{L}$ for all $n \in \mathbb{N}$. But $\mathfrak{F}[p^n] = \mathfrak{F}[\omega^{en}]$. Hence $\{z_1(en), \cdots, z_{h_{r,A}}(en)\}$ generates $\widehat{L}_{\mathfrak{F}}(p^n)$ over $\widehat{L}$. Thus using notation of previous section \[m_p(n) \leq h_{r,A} = \frac {h_A} {f}\].\\
\textbf{Lemma 3.2 :} i) Hypothesis as above. Let $d$ be dimension of $H$ as $\mathbb{Q}_p$ manifold. Let $d = \frac {h_Ah_p} {f}$ and $n_0 \in \mathbb{N}$ be such that $k \geq n_0$ implies $|I_k| = \text{vol}(H)p^{kd}$. Then $U_{n_0,A} \xhookrightarrow{} K_{n_0} \subseteq H$. \\ 
ii) Further assume that $n_0 \geq 2$. Then $\text{exp}(p^{n_0}[\omega]) \in K_{n_0} \subseteq H$. \\~\\
\textbf{Proof :} i)  From lemma-2.5(iii) $d \leq h_pm_p(n)$ for all $n \geq n_0 + 1$. Using the bound obtained before we conclude $\{z_1(en), \cdots z_{h_A/f}(en)\}$ is a generating set with smallest size for all $n \geq n_0 + 1$. Fix this choice of generating set.\\
Let $n \geq n_0 + 1$. Consider the map \[i_n : U_{n_0}/U_{n} \to G_{n,n_0} \] defined by $i_n(\tilde{a}) = \tilde{a}_{\sigma}$ where $\tilde{a}_{\sigma} \in G_{n,n_0}$ is determined by 
\[
\tilde{a}_{\sigma}(z_i(en)) = [a]z_i(en) \;\forall \,1 \leq i \leq h_A/f
\]
for $a \in U_{n_0,K}$ (here $\tilde{a}$ means class of $a$).\\
Such $\tilde{a}_{\sigma}$ exists by lemma-2.6$(i)$ and it is easy to verify that $i_n$ is a well-defined group homomorphism. Injectivity of $i_n$ follows from injectivity of $\Delta_{n,n_0}$. \\
Since $i_n$ is naturally defined one can pass to projective limit and obtain the desired inclusion. \\
ii) Proof is similar to proof of part $(i)$.\\
There is an endomorphism $[\omega] : \mathfrak{F}[p^k] \to \mathfrak{F}[p^k]$ for each $k \in \mathbb{N}$.\\
As before $\{z_1(en), \cdots z_{h_A/f}(en)\}$ is a minimal generating set for $\widehat{L}(p^n)$ over $\widehat{L}$.\\
A simple calculation shows \[\text{exp}(p^{n_0}[\omega]) = 1 + p^{n_0}\sum_{i \geq 1}a_i[\omega]^i\] where $a_i \in \mathbb{Z}_p$ for all $i \geq 1$ .\\
Thus by lemma-2.6(i) one can construct an element $\sigma \in G_{n,n_{0}}$ such that \[\sigma(z_i(en)) = \text{exp}(p^{n_0}[\omega])z_i(en)\] for all $1 \leq i \leq h_A/f$ and for all $n \geq n_0 + 1$. Passing to projective limit one obtains an element of $K_{n_0} \subseteq H$. $\square$ \\~\\
\textbf{Remark 3.3 :} i) Note that in proof of part $(i)$ we have constructed the map in such a way that the image of $U_{n_0,A}$ is indeed the image of $U_{n_0,A} \subseteq A$ under the inclusion map $A \to \text{End}_{\mathbb{Z}_p}(T_{p}(\mathfrak{F}))$. Here one should choose base to be $\{[\omega]^jz_{i} \,|\, 0 \leq j \leq e-1, 1 \leq i \leq \frac {h_A} {f} \}$ to get suitable matrix representation.   \\
ii) From lemma-3.2 it follows that if $d = \frac {h_ph_A} {f}$ then one has desired $A$-structure on $H$ and using the topological isomorphism one can put $A$-structure on $H_{\omega}$. \\~\\
We introduce a terminology :\\~\\
\textbf{Definition 3.4 :} Let $k$ be a complete discrete valuation field of characteristic $(0,p)$ and let $\mathfrak{F} \in \mathcal{F}_p(O_k)$. If $\mathfrak{F}$ is said to be \emph{good} over $k$ if $\text{End}_{O_k}(\mathfrak{F})$ is integrally closed in its fraction field.  \\~\\
\textbf{Remark 3.5 :} i) Let $\mathfrak{F} \in \mathcal{F}_{\pi}(O_{\widehat{L}})$. Then $\text{End}_{O_{\widehat{L}}}(\mathfrak{F}) = \text{End}(\mathfrak{F})$ and $\mathfrak{F}$ is good over $\widehat{L}$ (see \cite[Section-2]{some more results}).\\
iii) Let $\mathfrak{F} \in \mathcal{F}_{p}(O_{\widehat{L}})$. Then there is a $\mathfrak{G} \in \mathcal{F}_p(O_{\widehat{L}})$ such that $\text{End}_{O_{\widehat{L}}}(\mathfrak{F}, \mathfrak{G}) \neq 0$ and $\mathfrak{G}$ is good over $\widehat{L}$ (see \cite[3.2]{lubin2}). Thus their $p$-torsions generate the same fields and we only need to consider $\mathfrak{G}$ if we want to study the associated Galois group (see \cite[Remark-2.2]{some more results}).\\~\\
Consider the set-up in beginning of this section. Further let $A \xhookrightarrow{} \text{End}(\mathfrak{F})$ be onto. One can find such $A$ if and only if $\mathfrak{F}$ is good over $\widehat{L}$.\\ In this situation one can always show that $H_{\omega}$ has structure of $K(A)$ manifold and it is open in $\text{Gl}_{A}(T_{\omega}(A))$ by a result of Serre-Sen (see below).\\~\\  
One can generalize Proposition-2.4 to arbitrary local fields :\\~\\
\textbf{Proposition 3.6 :} Let $Y$ is a closed, smooth, analytic sub-manifold of $(A)^N$ and let $Y_n$ denote its image under reduction modulo $\omega^n$. Let $d_A$ be dimension of $Y$ and $\text{vol}(Y)$ be volume of $Y$ in the induced measure. Then $|Y_n| = \text{vol}(Y)p^{nd_Af}$ for large enough $n$ where $f$ is the degree of residue extension of $K(A)/\mathbb{Q}_p$.\\~\\
\textbf{Proof :} Similar to Serre's proof of Proposition-2.4 (see \cite[Section-3]{serre 3}). $\square$ \\~\\
With help of this proposition one can prove analogues of results in section-2. In particular if $H_{\omega}$ has structure of $K(A)$ manifold \[|\text{Gal}(\widehat{L}_{\mathfrak{F}}(\omega^{e(n+1)}) | \widehat{L}_{\mathfrak{F}}(\omega^{en}))| = p^{ed_Af}\] for large enough $n$. But $\mathfrak{F}[\omega^{en}] = \mathfrak{F}[p^n]$. Comparing with result in last section one obtains \[ed_Af = d\] where $d$ is dimension of $H_{\omega}$ as $\mathbb{Q}_p$ manifold and $d_A$ is dimension as $K(A)$ manifold.\\~\\
Hypothesis be as before ie $\mathfrak{F}$ is good over $\widehat{L}$ and $A = \text{End}(\mathfrak{F})$. Let $\mathfrak{h}$ denote lie algebra of $H$ thought as a subspace of $\text{End}_{\mathbb{Q}_p}(V_p)$ where $V_p = T_{p}(\mathfrak{F}) \otimes_{\mathbb{Z}_p} \mathbb{Q}_p$. $V_p$ has a natural $K(A)$ module structure, it is a $K(A)$ module of dimension $h_{r,A}$ and $K(A)$ can be identified with a sub-algebra of $\text{End}_{\mathbb{Q}_p}(V_p)$. Further Tate's main theorem on $p$-divisible group implies $K(A) = \text{End}_{\mathfrak{h}}(V_p)$ (\cite[Proposition 7]{serre 2}).\\
We shall use notation $V_{\omega} = T_{\omega}(\mathfrak{F}) \otimes_{A} K(A)$. \\~\\
Following results are due to Serre and Sen (\cite{sen2}) :\\~\\
\textbf{Lemma 3.7 :} $K(A)\mathfrak{h} \subseteq \mathfrak{h}$.\\~\\
\textbf{Proof :} This is consequence of Sen's theory of Hodge-Tate representation. See proof of \cite[Theorem-3]{sen2}. $\square$ \\~\\
One uses this lemma and analogous techniques from proof of Proposition 2.1 to show :\\~\\
\textbf{Proposition 3.8 :} $\text{End}_{K(A)}(V) = \mathfrak{h}$.\\~\\
\textbf{Proof :} See \cite[Theorem-3]{sen2}.\\ This proof requires hypothesis $K(A) = \text{End}_{\mathfrak{h}}(V_p)$. $\square$ \\~\\
Results above implies : \\~\\
\textbf{Proposition 3.9 :} $H_{\omega}$ is a lie-group over $K(A)$ of dimension $h_{r,A}^{2}$ which is open in $\text{Gl}_{A}(T_{\omega}(\mathfrak{F}))$.\\~\\ 
\textbf{Proof :} By Lemma 3.8 and Proposition 3.9, $\mathfrak{h}$ has a $K(A)$ module structure and as $K(A)$-module its dimension is $h_{r,A}^{2}$. Hence $H$ is a lie group over $K(A)$ of dimension $h_{r,A}^{2}$. Using the topological isomorphism as before, $H_{\omega}$ can be endowed with lie group structure over $K(A)$ of dimension $h_{r,A}^{2}$. Now proposition follows from following lemma and the fact that $\text{Gl}_{A}(T_{\omega}(\mathfrak{F}))$ is open subgroup of $\text{Gl}_{K(A)}(V_{\omega})$ containing $H_{\omega}$. $\square$ \\~\\
The following lemma is basic result from $p$-adic analysis :\\~\\
\textbf{Lemma 3.10 :} Let $k$ be a local field and $G$ be a topological subgroup of $\text{Gl}_{n}(k)$ which can be endowed with analytic structure over $k$ of dimension $n^2$. Then $G$ is an open subgroup . \\~\\
\textbf{Proof :} It is enough to show a neighbourhood of identity in $G$ contains an open set of $\text{Gl}_n(k)$.\\
Consider $G$ and $\text{Gl}_n(k)$ as lie groups over $\mathbb{Q}_p$.  Inclusion map is analytic as homomorphism of lie groups over $\mathbb{Q}_p$ (\cite[Part II, Ch V, Sec 9]{serre 1}). Since both sides have same dimension and there is no kernel, inclusion map is \'etale at $\text{Id}$ and hence an open map. $\square$ \\~\\ 
\textbf{Remark 3.11 :} i) Content of proposition-3.9 seems to be well-known. Lie group structure of $H_{\omega}$ over $K(A)$ can be derived in alternate manner after obvious identifications (see proof of Proposition A.23 in appendix).\\
ii) $d_A = h_{r,A}^{2}$ and $d = \frac {h_ph_A} {f}$ if $\mathfrak{F}$ is good and $A = \text{End}(\mathfrak{F})$.\\
iii) Assume that $H_{\omega}$ has $K(A)$-structure. From here one can prove image is open using suitable formalism of $\pi$-divisible groups. This is developed in appendix.\\~\\ 
\section{Some implications}
This section is devoted to collect some observations and answer a few questions posed by the author in previous articles \cite{some results}, \cite{some more results} which mostly follow from proposition-3.9. \\~\\
\textbf{Remark 4.1 :} i) Let $\mathfrak{F}, \mathfrak{G} \in \mathcal{F}_{p}(O_{\widehat{L}})(h)$ be good and $A_1, A_2$ be their endomorphism rings respectively. From remark-3.11(ii) it follows that $\text{dim}_{\mathbb{Q}_p}(H(\mathfrak{F})) = \text{dim}_{\mathbb{Q}_p}(H(\mathfrak{G}))$ if and only if $[K(A_1):\mathbb{Q}_p] = [K(A_2) : \mathbb{Q}_p]$.\\
ii) Let $\mathfrak{F} \in \mathcal{F}_p(O_{\widehat{L}})(h)$ be good and let $A$ denote its ring of endomorphisms. From remark-3.11(ii) it follows that $m_p(n) = \frac {h_A} {f}$ for large enough $n$.\\
iii) Further, let $\omega$ be a generator of the maximal ideal of $A$. Use $m_{\omega}(n)$ to denote the smallest size of a subset of $\mathfrak{F}[\omega^n] - \mathfrak{F}[\omega^{n-1}]$ which generate $K_{\mathfrak{F}}(\omega^n)$ over $K$. Clearly, $m_p(n) = m_{\omega}(en)$ for all $n$.\\
From proposition-3.6 and remark-3.11(ii) one can deduce that $m_{\omega}(n) = \frac {h_A} {f}$ for large enough $n$.\\
iv) By proposition-3.9 $H_{\omega}$ contains an open subgroup of $A^{\times}\text{Id} \subseteq \text{Gl}_{A}(T_{\omega}(\mathfrak{F})).$  Note that, $H$ contains same elements from image of $A^{\times}$ in $\text{Gl}_{\mathbb{Z}_p}(T_p(\mathfrak{F}))$. \\
This conclusion can be alternately derived from lemma-3.2 and remark-3.11(ii).\\~\\
Now we shall consider unramified group laws. Notation be as in introduction.\\
Assume that $\mu_{p^h-1} \subseteq O_E$. Note that for $\mathfrak{F} \in \mathcal{F}_{\pi}(O_E)(h)$, $\text{End}(\mathfrak{F}) = \text{End}_{O_E}(\mathfrak{F})$ and this ring is a complete, integrally closed subring of $O_E$ containing $A_{\pi}$. For simplicity, it will be denoted by $A_{\mathfrak{F}}$. Put $[k(A_{\mathfrak{F}}) : \mathbb{F}_p] = f_{\mathfrak{F}}$. Clearly, $f_{\pi}\,|\,f_{\mathfrak{F}}$ and $f_{\mathfrak{F}}\,|\,h$. Further $\pi$ is a generator for maximal ideal of $A_{\mathfrak{F}}$. \\
We shall be considering representation $\widehat{\rho}_{\pi}$ of $\text{Gal}(C|\widehat{L})$ ie the inertia subgroup, on $\text{Gl}_{A_{\pi}}(T_{\pi}(\mathfrak{F}))$. Note that $T_{\pi}(\mathfrak{F})$ is a free $A_{\pi}$ module of rank $h_{r,\pi}$. Thus $\text{Gl}_{A_{\pi}}(T_{\pi}(\mathfrak{F}))$ can be identified with $\text{Gl}(h_{r,\pi},A_{\pi})$. Since $T_{\pi}(\mathfrak{F})$ has a $A_{\mathfrak{F}}$ module structure, image of $\widehat{\rho}$ (denoted $H_{\pi}(\mathfrak{F}$)) lies inside $\text{Gl}(h_{r,A_{\mathfrak{F}}}, A_{\mathfrak{F}}) \subseteq \text{Gl}(h_{r,\pi}, A_{\pi})$.\\~\\
\textbf{Remark 4.2 :} i) By proposition-3.9, $H_{\pi}(\mathfrak{F})$ is a lie group over $K(A_{\mathfrak{F}})$. Hence it has analytic structure over $K_{\pi}$. This answers a question in remark-3.1.2 of \cite{some more results}. \\
ii) Further, as $K(A_{\mathfrak{F}})$ manifold $H_{\pi}(\mathfrak{F})$  has dimension $h_{r,A_{\mathfrak{F}}}^{2}$. Hence its dimension as $K_{\pi}$ manifold is \[h_{r,A_{\mathfrak{F}}}^2 \times \frac {f_{\mathfrak{F}}} {f_{\pi}} = \frac {h^2} {f_{\mathfrak{F}}f_{\pi}}.\]
This answers a question in remark-4.2.1 of \cite{some results}.\\    
iii) Let $\mathfrak{F}, \mathfrak{G} \in \mathcal{F}_{\pi}(O_E)(h)$. Clearly \[\begin{split} A_{\mathfrak{F}} = A_{\mathfrak{G}} \iff K(A_{\mathfrak{F}}) = K(A_{\mathfrak{G}})  \iff  \\ [K(A_{\mathfrak{F}}) : \mathbb{Q}_p] = [K(A_{\mathfrak{G}}) : \mathbb{Q}_p]  \iff f_{\mathfrak{F}} = f_{\mathfrak{G}}. \end{split} \] 
Thus remark-4.1(i) can be rephrased as \[\begin{split}\text{dim}_{\mathbb{Q}_p}(H(\mathfrak{F})) = \text{dim}_{\mathbb{Q}_p}(H(\mathfrak{G})) \iff \text{dim}_{K_{\pi}}(H_{\pi}(\mathfrak{F})) = \text{dim}_{K_{\pi}}((H_{\pi}(\mathfrak{G}))\\ \iff A_{\mathfrak{F}} = A_{\mathfrak{G}}. \end{split} \] 
iv) Note that $\text{dim}_{K_{\pi}}(H_{\pi}(\mathfrak{F})) = h_{r,\pi} \times \frac {h} {f_{\mathfrak{F}}}$ ie an integral multiple of $h_{r,\pi}$. Using remark-4.1(iii) one sees that if $m_{\pi}(n) = 1$ for large enough $n$, then $h = f_{\mathfrak{F}}$. \\ 
v) $H_{\pi}(\mathfrak{F})$ contains an open subgroup of $A_{\mathfrak{F}}^{\times}\text{Id} \subseteq \text{Gl}(h_{r,A_{\mathfrak{F}}}, A_{\mathfrak{F}}) $.\\~\\
Now we have an improved version a theorem from \cite{some results} (remark 1.1(iii)) : \\~\\
\textbf{Theorem 4.3 :} Let $\mathfrak{F} \in \mathcal{F}_{\pi}(O_E)(h)$ be an unramified group law. Further, assume that $\mu_{p^h - 1} \subseteq O_E$. Then, following are equivalent :\\
i) $\text{Gal}(E_{\mathfrak{F}}(\pi^{\infty})|E)$ is abelian,\\
ii) $m_{\pi}(n) = 1$ for large enough $n$,\\
iii) $\mathfrak{F}$ has $A_{\pi}^{(h_{r,\pi})}$ module structure extending $A_{\pi}$ module structure.\\~\\
\textbf{Proof :} Follows from theorem-2.3 of \cite{some results} and remark-4.2(iv). $\square$  
       
\section{Concluding remarks}
We conclude with some questions and observations about overall situation :\\~\\
\textbf{Remarks 5.1 :}
i) Let  $\mathfrak{F}, \mathfrak{G} \in \mathcal{F}_{p}(O_{\widehat{L}})(h)$. It is interesting to know when $\widehat{D}_p(\mathfrak{F}) = \widehat{D}_p(\mathfrak{G})$. Note that if they are in same isogeny class, this is the case. One can ask if the converse is true, ie. if $\widehat{D}_p(\mathfrak{F}) = \widehat{D}_p(\mathfrak{G})$ is $\mathfrak{F}$ and $\mathfrak{G}$ are isogenous ?\\
ii) One way to get information about containment like $\widehat{D}_p(\mathfrak{F}) \subseteq \widehat{D}_p(\mathfrak{G})$ is to compute the Galois group $\text{Gal}(C_p/\widehat{D}_p(\mathfrak{F}))$. Since $\widehat{D}_p(\mathfrak{F})/\widehat{L}$ is a lie extension, it is \emph{arithmetically pro-finite} and one can use `field of norm' machinery to compute the associated Galois group (see \cite{winten}). \\
In particular, if $\mathfrak{F} = \mathbb{G}_m$, there is an well-known formula for $\text{Gal}(\overline{\mathbb{Q}}_p|\mathbb{Q}_p(p^{\infty}))$ in terms of field of norms. \\
iii) Coleman has studied norm compatible systems in Lubin-Tate case using interpolation techniques (see \cite{col}). One would like to generalize these results for any unramified group law. Note that in this case one expects $h_{r, \pi}$ many power series in one variable whose combination will represent the element.\\~\\

\section{Acknowledgments}
I am thankful to Prof.C. Kaiser (Max Planck Institute for Mathematics, Bonn) for pointing out the result of proposition-2.3 in Serre's paper (see section-2) and Prof. L. Berger (ENS, Lyon) for pointing out the result of proposition-3.8 in Sen's paper (see section-3).\\  
I am grateful to Max Planck Institute for Mathematics, Bonn for hospitality and financial support.\\~\\  
\appendix
\begin{center}
\LARGE \textbf{Appendix} 
\end{center}
\vspace{1cm}
\section{Introduction}
The purpose of this appendix is to generalize the results of \cite{serre 2} for lie groups over arbitrary local fields. Since we are interested in representations arising from formal $A$-modules one should improve the formalism of Tate's work on $p$-divisible groups (see \cite{tate}) to prove a generalized `Hodge-Tate decomposition' for the Tate module which is a key ingredient for Serre's work.\\
We shall use notion of formal module schemes and follow Faltings work \cite{faltings} to prove a general version of classical Cartier duality which shall lead a generalized `Hodge-Tate decomposition'.  The results in \cite{serre 2} are easily generalized once we have a general `Hodge-Tate decomposition'.\\
Concept of $\mathcal{O}$-module schemes goes back work to Drinfeld and Faltings' paper (\cite{faltings}) does an elaborate treatment of the topic.
\section{Definition and first properties}
Let $\mathcal{O}$ be ring of integers of a local field and $\pi$ be a generator of the maximal ideal $\mathfrak{p}$ of $\mathcal{O}$. Use $K(\mathcal{O})$ to denote the fraction field of $\mathcal{O}$, $e$ and $f$ be ramification index and degree of residue extension associated with the extension with $K(\mathcal{O})/\mathbb{Q}_p$.\\~\\
\textbf{Definition A.1 :} i) Let $R$ be a commutative ring and let $G$ be a finite, flat, affine commutative group scheme over $R$. $G$ is said to be a \emph{finite $\mathcal{O}$-module scheme} if there is a homomorphism of rings \[ \mathcal{O} \to \text{End}_{R}(G,G).\]
ii) Further, assume that $R$ is an $\mathcal{O}$ algebra. A finite $\mathcal{O}$-module scheme is said to be \emph{strict} if $\mathcal{O}$ acts on $\text{Lie}(G)$ by scalar multiplication.\\~\\ 
\textbf{Remark A.2 :} i) In what follows we shall always assume $R$ is a $\mathcal{O}$ algebra and all $\mathcal{O}$ module schemes are strict. \\
ii) Note that $G = \text{Spec}(A)$ where $A$ is a locally free $R$-algebra of finite rank. Thus $\text{End}_R(G,G) = \text{End}_{R-\text{alg}}(A,A)$ as a set. Further if $f, g \in \text{End}_R(G,G)$ and $f^{*}$ and $g^{*}$ are corresponding elements in $\text{End}_{R-\text{alg}}(A,A)$ then $(f+g)*$ corresponds to $m \circ (f^{*} \otimes g^{*}) \circ \mu$ where $\mu$ and $m$ are co-multiplication and algebra multiplication maps and $(f \circ g)^{*}$ corresponds to $g^{*} \circ f^{*}$.\\
iii) Homomorphism of two $\mathcal{O}$-modules are defined as homomorphism of group schemes which commute with module structure. If $G$ and $H$ are two finite $\mathcal{O}$ module schemes then we shall use the notation $\text{Hom}^{\mathcal{O}}_{R}(G,H)$ to describe the collection of all homomorphism of $\mathcal{O}$-module schemes.\\
iv) Using functor point of view for schemes the definition above implies that a finite $\mathcal{O}$-module schemes is a representable functor from the category of $R$ algebras to the category $\mathcal{O}$-modules with finitely many elements.  \\~\\
Let $R$ be a noetherian, local, complete $\mathcal{O}$ algebra and $G$ be a finite, flat, affine, commutative group scheme over $R$. Then there is an exact sequence \[ 0 \to G^{0} \to G \to G^{\text{et}} \to 0 \tag{A.1}\] of group schemes (see \cite[1.4]{tate}), where $G^{0}$ is connected component of identity which is closed, flat subgroup of $G$ and quotient $G^{\text{et}}$ is an \'etale, finite group scheme over $R$.\\ If $G$ has structure of $\mathcal{O}$ module it is easy to check that $G^{0}$ and $G^{\text{et}}$ inherits structure of $\mathcal{O}$ module. Note that, $\text{Lie}(G) = \text{Lie}(G^0)$ and $\text{Lie}(G^{\text{et}}) = 0$. Hence if $G$ has strict module structure (which is always the case below), so do $G^0$ and $G^{\text{et}}$.\\~\\
If $G$ is a finite, flat group scheme over an $\mathcal{O}$ algebra $R$ which has connected spectrum, we use $|G|$ to denote rank of locally free sheaf defining $G$. \\~\\  
\textbf{Definition A.3 :} Let $h \in \mathbb{N}$ and $R$ be a connected $\mathcal{O}$ algebra. A $\pi$-divisible group over $R$ is defined as an inductive system $(G_{\nu}, i_{\nu})_{\nu \geq 0}$ of finite $\mathcal{O}$-module schemes over $R$ satisfying :\\
i) $|G_{\nu}| = p^{\nu h}$ for each $\nu \geq 0$ for some fixed integer $h$,\\
ii) For each $\nu \geq 0$, there is an exact sequence \[0 \to G_{\nu} \xrightarrow{i_{\nu}} G_{\nu + 1} \xrightarrow{[\pi^{\nu}]} G_{\nu}.\] \\
\textbf{Remark A.4 :} i) $h$ is said to be height of $G$.\\
ii) A homomorphism of two $\pi$-divisible groups $(G_{\nu}, i_{\nu})$ and $(H_{\nu}, j_{\nu})$ is a sequence $f = (f_{\nu})$ such that $f_{\nu} : G_{\nu} \to H_{\nu}$ is a morphism of $\mathcal{O}$-module schemes over $R$ such that $ f_{\nu +1}i_{\nu} = j_{\nu+1}f_{\nu}$ for all $\nu$. \\
iii) Note that $i_{\nu} : G_{\nu} \to G_{\nu+1}$ identifies $G_{\nu}$ with a closed sub-scheme of $G_{\nu+1}$. Iterating we have a closed immersion $i_{\nu, \mu} : G_{\nu} \to G_{\nu + \mu}$. Further image of this map coincides with kernel of $[\pi^{\nu}] : G_{\nu + \mu} \to G_{\nu + \mu}$. Image of $[\pi^{\nu}]$ can be identified with $G_{\mu}$. Considering orders one has an exact sequence \[0 \to G_{\mu} \xrightarrow{i_{\mu,\nu}} G_{\mu+\nu} \xrightarrow{[\pi^{\mu}]} G_{\nu} \to 0.\]\\
iv) Consider the inductive system $(G_{e\nu}, i_{\nu}(p))_{\nu \geq 0}$ where \[i_{\nu}(p) = i_{e(\nu +1) - 1} \circ \cdots \circ i_{e\nu}.\] It is easy to check that this defines a $p$-divisible group $G_{p} = \varinjlim G_{e\nu}$ of height $eh$. Note that $\mathcal{O} \xrightarrow{} \text{End}_{R}(G_{p})$. Thus one can consider $[\pi]$ torsion sub-modules and recover the associated $\pi$-divisible group.\\
v) Put $|\mathcal{O}/\pi\mathcal{O}| = p^f$. In analogy of formal groups we would like to have $f | h$. If $L$ is a separably closed field in $\text{Affine}_R$ of characteristic $\neq p$, $|G_{\nu}| = |G_{\nu}(L)|$. A cardinality argument shows in this case $f \,|\,h$ and as $\mathcal{O}$ module $G_{\nu}(L) \cong (\mathcal{O}/\pi^\nu \mathcal{O})^{h_r}$ where $h_r = \frac {h} {f}$. Such a field exists if $R$ is a domain of char $0$ which is always the case below. In this situation $h_r$ is called $\pi$-height.\\~\\
For rest of the section $R$ be a noetherian, local, complete domain containing $\mathcal{O}$ of char $(0,p)$.\\~\\
From connected-\'etale exact sequence $(A.1)$ one concludes that $(G^{0}_{\nu},i_{\nu})_{\nu \geq 0}$ forms a $\pi$-divisible groups. We shall write this group as $G^{0}$.\\
Similarly $(G_{\nu}^{\text{et}},i_{\nu})_{\nu \geq 0}$ is also a $\pi$-divisible groups denoted $G^{\text{et}}$ and one has exact sequence \[0 \to G^{0} \to G \to G^{\text{et}} \to 0.\]
$G$ is said to be connected if $G = G^{0}$ and \'etale if $G = G^{\text{et}}$.\\~\\
Next we shall show that category of connected $\pi$-divisible groups over $R$ is equivalent to category of formal $\mathcal{O}$-modules over $R$ with property that $[\pi]$ defines an isogeny. \\
Let $\mathfrak{F}$ be a formal $\mathcal{O}$ module over $R$ of dimension $n$. Put $A = R[[\underline{X}]]$ where $\underline{X}$ is a short-hand for $X_1, \cdots, X_n$. Using the module structure one obtains a $R$ algebra morphism $\phi : A \to A$ defined by $X_i \to [\pi]_i(\underline{X})$. This turns $A$ into a $A$ module. \\~\\
\textbf{Definition A.5 :} $\mathfrak{F}$ is said to be \emph{divisible} if $\phi$ turns $A$ into a finite free $A$ module. \\~\\
Put $J_{\nu} = \phi^{\nu}((\underline{X}))A$ and $A_{\nu} = A/J_{\nu}$ for $\nu \geq 0$. Define $G_{\nu}(\mathfrak{F}) = \text{Spec}(A_{\nu})$. Clearly $G_{\nu}(\mathfrak{F})$ inherits $\mathcal{O}$ module structure from $\mathfrak{F}$ and one has a $\pi$-divisible  group $G(\mathfrak{F}) = \varinjlim G_{\nu}(\mathfrak{F})$. Note that $f$ times $\pi$-height of $G(\mathfrak{F})$ is rank of $A$ as $A$ module. One can check $G_{\nu}(\mathfrak{F})$ is connected for $\nu \geq 0$ and thus we have a connected $\pi$-divisible group over $R$.\\~\\
Conversely, if one has a connected $\pi$-divisible group $G$ consider the associated $p$-divisible group $G_p$. By Tate's construction there is a divisible formal group law $\mathfrak{F}$ over $R$ such that $G_p$ arises from $\mathfrak{F}$ in the sense described above. Note that $\mathcal{O} \xhookrightarrow{} \text{End}_R(\mathfrak{F}, \mathfrak{F})$ and acts on $\text{Lie}(\mathfrak{F})$ (constructed from formal group perspective) by scalars. Hence we have a formal $\mathcal{O}$ module over $R$. Since $\mathfrak{F}$ is divisible as $\mathbb{Z}_p$ module it is also divisible as $\mathcal{O}$ module. \\~\\
\textbf{Proposition A.6 :} The correspondence $\mathfrak{F} \to G(\mathfrak{F})$ is an equivalence of categories between divisible formal $\mathcal{O}$ modules and connected $\pi$-divisible groups over $R$.\\~\\
\textbf{Proof :} Follows from discussion above. See also \cite[Section-2.2]{tate}. $\square$  \\~\\
\textbf{Remark A.7 :} i) With terminology as in introduction, any unramified group law in one variable is divisible (\cite[Section 2]{some results}). In particular one can consider a Lubin-Tate group law defined over $\mathcal{O} \subseteq R$. It is divisible over $\mathcal{O}$ and hence over $R$. We shall denote the associated $\pi$-divisible group as $G_{\text{LT}}$.\\
ii) Proposition A.6 gives natural examples of $\pi$-divisible groups. As consequence, any $\mathcal{O}$ module of finite height in one variable is divisible over $R$.\\
iii) Let $G$ be a $\pi$-divisible group over $R$. Dimension of $G$ is defined as dimension of formal group associated to $G^{0}$.\\~\\
Let $\mathfrak{F}$ and $\mathfrak{G}$ be formal $\mathcal{O}$-modules over $R$. Since characteristic of $R = 0$ it can be shown, $\text{Hom}_R(\mathfrak{F}, \mathfrak{G}) = \text{Hom}_R^{\mathcal{O}}(\mathfrak{F}, \mathfrak{G})$ (\cite[21.1.4]{haz}). \\~\\
\textbf{Lemma A.8 :} Let $G$ and $H$ be two $\pi$-divisible groups and $G_p$, $H_p$ be corresponding $p$-divisible groups. If $H$ is connected, \[ \text{Hom}_{R}(G_p, H_p) = \text{Hom}_R^{\mathcal{O}}(G_p, H_p). \] 
\textbf{Proof :} It is enough to show  that any homomorphism of $p$-divisible groups $f : G_{p} \to H_p$ is a homomorphism of $\mathcal{O}$ modules. \\
For $a \in \mathcal{O}$, consider $f_a : G_p \to H_p$ defined by $f_a = [a]\circ f - f\circ[a]$. Clearly $f_a$ is a homomorphism of groups. Let $G_p^{0}$ be the connected component  of $G_p$. Using the result above, $f_a$ is identically zero when restricted to $G_p^{0}$. Thus we have a homomorphism of groups $G_{p}^{\text{et}} \to H_p$. By a result of Waterhouse \cite[Lemma-2.2]{water} such a homomorphism must be trivial since $H_p$ is connected. So $f_a = 0 $ for each $a \in \mathcal{O}$. \\ 
This proves the lemma. $\square$ \\~\\
The last topic of this section is tangent space of $\pi$-divisible groups. For this part we shall assume $R$ is a DVR with perfect residue field.\\
In what follows $L$ is completion of an algebraic extension of $K(R)$, the fraction field of $R$ and $O_L$ is ring of integers of $L$ and $\mathfrak{m}_L$ be its maximal ideal.\\
Let $G$ be a $\pi$-divisible group over $R$. We use notations $G(O_L)$, $\Phi_{\pi}(G)$ and $T_{\pi}(G)$ in sense of \cite[Section 2.4]{tate}. If $L$ is algebraically closed $\Phi_{\pi}(G)$ can be identified with group of torsion points of $G(O_L)$. \\~\\
Let $G$ be a $\pi$-divisible group and $G^{0}$ be its connected component of identity. Let $A_0$ be the co-ordinate ring for $G^0$. One can define the tangent space $t_G$ in terms of left invariant derivations of $A_0$. $S$ valued points of $t_{G}$ are given by $R$-linear maps $\tau : A_0 \to S$ such that $\tau(fg) = f(0)\tau(g) + \tau(f)g(0)$ for any $R$ algebra $S$. Note that $A_0 = R[[X_1, \cdots, X_n]]$ for some $n \in \mathbb{N}$. Put $I_0 = (X_1, \cdots, X_n)$. Then \[ t_{G}(S) = \{ f : I_0/I_0^{2} \to S \,|\, f \,\text{is}\,R\,\text{linear} \}. \]
It is a free module of rank $n$ over $S$ and this rank equals dimension of $G^{0}$.\\
For all $a \in \mathcal{O}$, one has map $[a] : A \to A$ defined by $X_i \to a_i(X)$. Since $a_i(X) = aX_i + \text{higher degree terms},$ $t_{G}(S)$ has a natural $\mathcal{O}$ module structure.\\
Consider the case $S = L$. One has logarithm map $\text{log} : G(O_L) \to t_{G}(L)$ given by \[ \text{log}(x)(f) = \lim_{i \to \infty} \frac {f\circ [\pi^i](x) - f(0)} {\pi^i}\]
where $x \in G(O_L)$ and $f \in A_0$.\\
Note that for any $x \in G(O_L)$, $[\pi]^{i}(x) \in G^{0}(O_L)$ for large enough $i$. So one can use the power series logarithm map coming from formal group and define the logarithm map. Using this point of view or otherwise it is easy to see that $\text{log}$ is  $\mathcal{O}$-linear with respect to $\mathcal{O}$ module structure on both sides. Further, log induces a local isomorphism between the groups on both sides. Thus we have an isomorphism of $\mathcal{O}$ modules \[ G(O_L) \otimes_{\mathcal{O}} K(\mathcal{O}) \xrightarrow{ \cong} t_{G}(L).\] \\
\textbf{Remark A.9 :} If valuation on $L$ is discrete then image of $\text{log}$ is a $O_L$ lattice in $t_{G}(L)$. If $L$ is algebraically closed then image is onto.

\section{Duality for $\pi$-divisible groups }
The goal of this section is to describe a dual of $\pi$-divisible groups which behave well with respect to action of endomorphism ring. This duality is similar to classical Cartier duality and quite important for rest of discussion.\\~\\
Let $R$ be an $\mathcal{O}$-algebra and $G$ be a affine, finite flat group scheme over $R$. In \cite[Section 2]{faltings} introduces an auxiliary construction and defines notion of strict $\mathcal{O}$ action which is stronger condition than the notion introduced earlier. Further, he defines strict homomorphism between two group schemes. We need this notion to construct duals.\\
However, kernels of multiplication by powers of $[\pi]$ in divisible formal $\mathcal{O}$ modules are finite, flat group schemes with strict $\mathcal{O}$ action in this sense. (\cite[Section 3]{faltings}) \\~\\     
\textbf{Construction A.10} \emph{(Generalized Cartier duality)} (\cite[Section-5]{faltings})\\
i) Let $G$ be a finite $\mathcal{O}$ module scheme over noetherian ring $R$ which has strict $\mathcal{O}$ action. $G_{LT}$ be Lubin-Tate module over $R$ obtained by base change. Since $G$ is finite there is an integer $h$ such that $\pi^h$ annihilates $G$. Clearly an $\mathcal{O}$-module homomorphism from $G$ into $G_{LT}$ factors through $G_{LT}[\pi^h]$.\\
ii) Define a functor from category of $R$ algebras to sets by $S \to \text{Hom}_{\text{strict}}(G_S, G_{LT,S})$ (ie the set of strict $\mathcal{O}$-module homomorphisms). This functor is representable by a finite $\mathcal{O}$ module-scheme $G'$ with strict $\mathcal{O}$ action which is defined over $R$ and of same rank as $G$. \\
iii) $(G')'$ isomorphic to $G$ and from definition it is clear that the association $G \to G'$ is functorial.\\
iv) Further $R$ be connected and $G$ be a $\pi$-divisible group of height $h$ such that each $G_{\nu}$ has strict $\mathcal{O}$ action. Such $\pi$-divisible group is said to have strict $\mathcal{O}$ action. Then one can construct another $\pi$-divisible group $G' = \{G'_{\nu}\}_{\nu \geq 0}$ defined over $R$ dualizing multiplication by $[\pi]$ map. It is easy to see that $G'$ has height $h$ and it has strict $\mathcal{O}$ action. \\~\\
Following proposition is analogue of a formula due to Tate :\\~\\
\textbf{Proposition A.11 :} $R$ be as in previous section and $G$ be a $\pi$-divisible group of height $h$ defined over $R$ with strict $\mathcal{O}$ action.\\
Then \[\text{dim}(G) + \text{dim}(G') = h_r \] where $h_r = \frac {h} {f}$, $\pi$-height of $G$.\\~\\
\textbf{Proof :} Write $q = p^f$. Let $A$ be a noetherian $\mathcal{O}$ algebra which is annihilated by $\pi$. For any $A$-scheme $X$ let $X^{(q)}$ denote its base change by $\text{Frob}_q$, the $q$-th power map on $A$. Then one has relative Frobenius $\text{Frob}_{\pi} : X \to X^{(q)}$ .\\
If $H$ is a finite $\mathcal{O}$ module scheme with strict $\mathcal{O}$ action, there exists a transfer $\text{Ver}_{\pi} : H^{(q)} \to H$ such that $\text{Frob}_{\pi}\circ \text{Ver}_{\pi} = [\pi]$ and $\text{Ver}_{\pi} \circ \text{Frob}_{\pi} = [\pi]$.\\
This $\text{Ver}_{\pi}$ is defined as adjoint of $\text{Frob}_{\pi}$ via duality. (\cite[Section 7]{faltings}).\\
Now $A$ be connected and $H = (H_{\nu},i_{\nu})$ be a $\pi$-divisible group of height $h$ with strict $\mathcal{O}$ action. One can construct another $\pi$-divisible group $H^{(q)} = (H_{\nu}^{(q)},i_{\nu}^{(q)})$ of height $h$ with strict $\mathcal{O}$ action and two isogenies $\text{Frob}_{\pi} : H \to H^{(q)}$ and $\text{Ver}_{\pi} : H^{(q)} \to H$ such that compositions are multiplication by $[\pi]$. One has an exact sequence \[0 \to \text{Ker}(\text{Frob}_{\pi}) \to \text{Ker}([\pi]_{H}) \to \text{Ker}(\text{Ver}_{\pi}) \to 0.\]     
Let $A$ be residue field of $R$ and think $G$ and $G'$ as $\pi$-divisible groups over $A$. Note that $\text{Ker}([\pi]_{G}) = G_1$, $|\text{Ker}(\text{Frob}_{\pi,G})| = q^{n}$, $|\text{Ker}(\text{Ver}_{\pi,G}|_{G_1^{(q)}})| = |\text{Coker}(\text{Ver}_{\pi,G}|_{G_1^{(q)}})|$ and later group is dual to $\text{Ker}(\text{Frob}_{\pi,G'})$ which has rank $q^{n'}$.\\
The identity follows from exact sequence considering rank. $\square$ 
\section{Hodge-Tate decomposition} 
Let $T$ be a free $\mathcal{O}$ module of finite rank and assume that $R$ is a complete $DVR$ containing $\mathcal{O}$ of char $(0,p)$ with perfect residue field. Put $V(K(\mathcal{O})) = T \otimes_{\mathcal{O}} K(\mathcal{O})$. Let $C$ denote completion of $\overline{K(R)}$, an algebraic closure of $K(R)$ and let $G_{K(R)}$ be the absolute Galois group of $K(R)$ and $\chi : G_{K(R)} \to K(R)^{\times}$ be continuous homomorphism.\\~\\
One can define twisted action of $G_{K(R)}$ on $C$ by $g_{\chi}(x) = \chi(g)g(x)$. We shall denote this $G_{K(R)}$ module as $C(\chi)$. First, we need some preliminary results which generalize Tate's computations of twisted cohomology. \\~\\
\textbf{Definition A.12 :} Let $\chi : G_{K(R)} \to K(R)^{\times}$ be a continuous character and $K_{\chi} = {\overline{K(R)}}^{\text{ker}(\chi)}$ be fixed field of $\text{ker}(\chi)$. $\chi$ is said to be \emph{Sen character} if there is a finite extension $K_{0}/K(R)$ in $K_{\chi}$ such that $K_{\chi}/K_0$ is a totally ramified pro-$p$ extension and its Galois group is identified by $\chi$ with $p$-adic lie group (over $\mathbb{Q}_p$) of dimension $\geq 1$ contained in $K(R)^{\times}$. \\~\\   
\textbf{Example A.13 :} Let $L = \bigcup_{n\geq 1}K(R)(G_{LT}[\pi^n])$. The character induced by Galois group of $L/K(R)$ is a Sen character. We normalize it as follows :\\
Let $z \in T_{\pi}(G_{LT})$ be a generator of $\mathcal{O}$ module. $\chi$ is defined by $g.z = \chi(g).z$ for all $g \in G_{K(R)}$. $\chi$ has image in $\mathcal{O}^{\times}$ and image is infinite. \\~\\
Let $K_0/K(R)$ be a finite extension and $\chi$ be a Lubin-Tate character over $K_0$. One can estimate valuation of discriminant, define normalized trace map and compute twisted cohomology (\cite[Section - 3]{tate}, also \cite[A.4.2]{fon} for stream-lined version of treatment). Only thing to note is results on totally ramified $\mathbb{Z}_p$ extension generalize to Lubin-Tate extensions since its Galois group contains an open sub-group isomorphic to $\mathcal{O}$ corresponding to totally ramified extension and tower of fields has good ramification theoretic description.\\~\\
\textbf {Proposition A.14 :} Let $\chi$ be Lubin-Tate character. \\
i) $H^{1}_{\text{cont}}(G_{K(R)}, C(\chi)) = 0$. \\
ii) $H^{0}_{\text{cont}}(G_{K(R)}, C(\chi)) = 0$. \\~\\
\textbf{Proof (Sketch) :} Let $X$ be the completion of $K_{\chi}$. Following Tate one can show $H^{1}_{\text{cont}}(G_{K(R)}, X(\chi)) = H^{0}_{\text{cont}}(G_{K(R)}, X(\chi)) = 0$. Now the proof follows from approximation of continuous cochains and general results in group cohomology (\cite[Section 3]{tate}).  $\square$\\~\\
\textbf{Remark A.15 :} i) Lubin-Tate group can be replaced by any $\mathcal{O}$ module with endomorphism ring of full height defined over ring of integers of a local field contained in $R$.\\
ii) In proposition A.14, $\chi$ can be replaced by $\chi^{-1}$.\\
iii) Does lemma A.14 holds for any Sen character ?\\~\\
Let $\chi : G_{K(R)} \to K(R)^{\times}$ be a continuous character such that $\chi^{i}$ is nontrivial for all $i \in \mathbb{Z}$  and \[V(C) = V(K(\mathcal{O})) \otimes_{K(\mathcal{O})} C.\] 
For integers $i \in \mathbb{Z}$ define 
\[\begin{split}
V_{i}(C) = \{ x \in V(C) \,|\, g.x = \chi(g)^{i}x \},\\
V(C)[i] = V_{i}(C) \otimes_{K(R)} C.
\end{split}\]    
Note that, $V_i(C)$ is $K(R)$-vector space and $V(C)[i]$ is $G_{K(R)}$ stable. The injection $V_{i}(C) \xhookrightarrow{} V(C)$ extends to $C$-linear injection $\epsilon_{i} : V(C)[i] \to V(C)$. This gives a $C$-linear map $ s : \bigoplus_{i \in \mathbb{Z}} V(C)[i] \to V(C)$. Following Serre one can show that $s$ is injective (see \cite[2.4]{serre 2}).\\~\\
\textbf{Definition A.16 :} $V$ is said to be \emph{$\chi$-Hodge-Tate module} if $s$ is onto ie $V(C) = \bigoplus_{i\in \mathbb{Z}} V(C)[i]$. \\~\\
\textbf{Remark A.17 :} i) To have all $\chi^{i}$ nontrivial it is enough to ensure $\chi$ has infinite image.\\
ii) If $V$ is $\chi$-Hogde-Tate then its dual $V^{*}$ is also $\chi$-Hodge-Tate and $V^{*}(C)[i] \cong V(C)[-i]$ as $G_{K(R)}$ module.\\
 \\~\\
Let $G$ be a $\pi$-divisible group over $R$ with strict $\mathcal{O}$ action and $G'$ be its dual as constructed in last section. We shall show that the $\pi$-adic Tate module $T_{\pi}(G)$ thought as $\mathcal{O}$ module, has $\chi$-HT decomposition for a suitable $\chi$.\\
There is a non-degerate pairing of $\mathcal{O}$ modules \[G_{\nu}(C) \times G'_{\nu}(C) \to G_{LT,\nu}(C)\] for each $\nu \geq 0$. This paring is compatible with action of $G_{K(R)}$ on both sides and taking injective limit in one co-ordinate and projective limit on other co-ordinate ie it induces a $\mathcal{O}$ linear, Galois invariant isomorphism \[T_{\pi}(G)\cong \text{Hom}_{\mathcal{O}}(T_{\pi}(G'), T_{\pi}(G_{LT})).\]\\
\textbf{Proposition A.18 :} As $G_{K(R)}$-module and $C$-vector space \[ \text{Hom}_{\mathcal{O}}(T_{\pi}(G), C) = t_{G'}(C) \oplus t^{*}_{G}(C)\otimes_{C}\text{Hom}_{\mathcal{O}}(H, C)\]     
where $H = T_{\pi}(G_{LT})$ and $t^{*}_G(C)$ is co-tangent space of $G$ at identity. $G_{K(R)}$ acts on homomorphisms via $(\sigma f)(x) = \sigma (f(\sigma^{-1}x))$. \\~\\
\textbf{Proof :} Let $W = \text{Hom}_{\mathcal{O}}(T_{\pi}(G), C)$, $W' = \text{Hom}_{\mathcal{O}}(T_{\pi}(G'), C)$ and $Y = \text{Hom}_{\mathcal{O}}(H, C)$. \\
We have an isomorophism of $G_{K(R)}$ as well as $\mathcal{O}$ modules \[T_{\pi}(G') \xrightarrow{\cong} \text{Hom}_{\mathcal{O}}(T_{\pi}(G), H).\] 
This gives a non-degenerate bi-linear Galois invariant pairing \[W \times W' \to Y.\] 
Now consider maps 
\[\begin{split}
d\alpha : t_{G}(C) \to W' \\
d\alpha' : t_{G'}(C) \to W 
\end{split}\]
as in \cite[Section-4]{tate}. One can imitate Tate's proof of proposition-11, step-6 and show these maps are injective. Further these are morphism of $G_{K(R)}$ modules.\\
Put $X = t_{G'}(C)$ and $X' = \text{Hom}_{C}(t_{G}(C),Y)$. By the pairing above, $X'$ can be thought as a subspace of $W$. First, we would like to show $W = X \oplus X'$ as $C$-vector space.\\
Note that $Y = C(\chi^{-1})$ as $G_{K(R)}$ module, where $\chi : G_{K(R)} \to K(R)^{\times}$ is a Lubin-Tate character. In particular $Y^{G_{K(R)}} = 0$. Note that $t_{G}(K(R)) \subseteq W'^{G_{K(R)}}$ and $t_{G'}(K(R)) \subseteq W^{G_{K(R)}}$. The bilinear paring above pairs $W^{G_{K(R)}} \times W'^{G_{K(R)}}$ into $Y^{G_{K(R)}}$. So $t_{G'}(K(R))$ and $t_{G}(K(R))$ are mutually orthogonal with respect to the pairing. Now from definition of $X'$ it follows that $X$ and $X'$ are linearly disjoint.\\
Note that $\text{dim}_{C}(W) = \frac {h} {f}$ where $h$ is height of $G$ and $f$ is residue degree of $\mathcal{O}$. \\
$\text{dim}_{C}(t_{G}(C)) = \text {dim of}\;G$ and $\text{dim}_{C}(t_{G'}(C)) = \text{dim of}\;G'$. From proposition A.11 $W = X \oplus X'$ as $C$-vector spaces.\\
Now one needs to show that the splitting above respects $G_{K(R)}$ module structure.\\
As in \cite[Corollary-2 of theorem-3]{tate}, it is enough to check \[H^{1}_{\text{cont}}(G_{K(R)}, C(\chi)) = 0\] which follows from proposition A.14. $\square$\\~\\
From the proposition one immediately deduces :\\~\\
\textbf{Corollary A.19 :} $T_{\pi}\otimes_{\mathcal{O}}K(\mathcal{O})$ is $\chi$-HT where $\chi$ is inverse of Lubin-Tate character.\\~\\
Here $V(C)[0] = t_{G'}(C)$ and $V(C)[1] = t^{*}_{G}(C)$ and these spaces have dimension equal to $\text{dim}(G')$ and $\text{dim}(G)$ respectively.
\section{Generalization of Serre's result} 
Let $V$ be a finite dimensional vector-space over $K(\mathcal{O})$ and let  \[\chi : G_{K(R)} \to K(R)^{\times}\] be a continuous character and\[\rho : G_{K(R)} \to \text{Gl}_{K(\mathcal{O})}(V)\] be a continuous representation such that :\\
i) Image of $\rho$ is a lie-group over $K(\mathcal{O})$,\\
ii) $\chi^i$ is nontrivial for all $i \in \mathbb{Z}$ and $V$ is $\chi$-HT with $V(C) = V(C)[0] \oplus V(C)[1]$,\\
iii) $V$ is a semisimple $G_{K(R)}$ module.\\~\\
Let $G = \text{Im}(G_{K(R)})$ and $\mathfrak{g}$ be corresponding lie-subalgebra of $\text{End}_{K(\mathcal{O})}(V)$. $G_{\text{alg}}$ be algebraic envelope of $G$ in $\text{Gl}_V$, thought as algebraic group over $K(\mathcal{O})$. By hypothesis $(iii)$ it follows that $V$ is a semi-simple $\mathfrak{g}$ module. $V$ is said to be absolutely simple as $\mathfrak{g}$ module if commutant of $\mathfrak{g}$ in $\text{End}_{K(\mathcal{O})}(V)$ is scalars.\\
Put $n_0 = \text{dim}_{C}(V(C)[0])$ and $n_1 = \text{dim}_{C}(V(C)[1])$.\\~\\ 
\textbf{Proposition A.20 :} i) Assume that $V$ is absolutely simple $\mathfrak{g}$ module and $n_0$ and $n_1$ are co-prime. Then $G_{\text{alg}} = {\text{Gl}}_{V}$.\\
ii) Let $\text{dim}_{K(\mathcal{O})}(V) = 1$. Then $G_{\text{alg}} = \text{Gl}_V$ if $G$ is infinite. \\~\\
Part $(ii)$ is clear. Proof of part $(i)$ is same as Serre's proof in \cite[Section-3 and Section-4]{serre 2}. We shall only mention few key points :\\~\\
\textbf{Claim A.21 :} Define $\phi : C^{*} \to \text{Aut}(V(C))$ where $\phi(\lambda) : V(C) \to V(C)$ is given by $\lambda^{i}$ on $V(C)[i]$. Put $\Phi = \text{Im}(\phi)$. Then $\Phi \subseteq G_{\text{alg}}(C)$. \\~\\
\textbf{Proof :} This is analogue of Serre's theorem-1 from section-3. Note that, the argument rests on existence of HT decomposition which easily goes over to $\chi$-HT modules. $\square$ \\~\\
Now one can apply proposition-5 from section-4 in \cite{serre 2} directly and conclude proposition A.20. \\~\\
We need a result about algebraic envelope :\\
Assumptions be as in beginning of this section. Since $V$ is a semisimple $G$ module, one can write $\mathfrak{g} = \mathfrak{c} \times \mathfrak {s}$ where $\mathfrak{c}$ is center and $\mathfrak{s} = [\mathfrak{g}, \mathfrak{g}]$ is semisimple.\\~\\
\textbf{Lemma A.22 :} i) If $\mathfrak{c}$ is algebraic lie algebra, $G_{\text{alg}}$ is equal to $\mathfrak{g}$ and $G$ is an open subgroup of $G_{\text{alg}}(K(\mathcal{O})).$\\
ii) If $\mathfrak{g}$ is absolutely simple, $\mathfrak{c}$ is algebraic.\\~\\
\textbf{Proof :} Part $(i)$ is proved in (\cite[Section 1]{serre 2}). Note that here base field is $K(\mathcal{O})$.\\
Part $(ii)$ follows from the fact that $\mathfrak{c}$, under assumption is either $0$ or consists of scalars both of which are algebraic.\\~\\       
\textbf{Proposition A.23 :} Notation be as in section-3. Let $\mathfrak{F} \in \mathcal{F}_{p}(O_{\widehat{L}})(h)$ be good with endomorphism ring $A$. Consider the associated Galois representation \[\rho : G_{\widehat{L}} \to \text{Gl}_{h_r}(K(A))\] where $h_r = \frac {h} {f}$, underlying vector space is $T_{\omega}(\mathfrak{F}) \otimes_{A} K(A)$ and we are fixing a base for $T_{\omega}(\mathfrak{F})$. \\
Then $\text{Im}(\rho)$ is an open sub-group of $\text{Gl}_{h_r}(K(A))$.\\~\\
\textbf{Proof :} First, we shall verify condition $(i)$ above. Let $\{z_1,\cdots, z_{h_r}\}$ be corresponding fixed base. Think these as elements of $T_{p}(\mathfrak{F})$. Clearly it is a $A$ base for $T_p(\mathfrak{F})$. Further we have an identification $\text{End}_{\mathbb{Z}_p}(T_p(\mathfrak{F}))$ with $\text{End}_{h_r}(A)$. One can identify image $H$ of Galois representation from section 3 with image $G$ above as endomorphisms since Galois group acts $A$ linearly. Let $\mathfrak{h} \subseteq \text{End}_{\mathbb{Q}_p}(V_p)$ be lie algebra corresponding to $H$. By lemma 3.7 it is closed under $K(A) = A \otimes_{\mathbb{Z}_p} \mathbb{Q}_p$. Hence corresponding subgroup in $\text{End}_{h_r}(K(A))$ is also a lie sub-algebra over $K(A)$. Looking at representation of this lie sub-algebra on underlying vector space one concludes that $G$ contains an analytic open subgroup whose lie algebra is this lie subalgebra. Thus $G$ has structure of analytic subgroup of $\text{Gl}_{h_r}(K(A))$ over $K(A)$.\\
Condition $(ii)$ follows from corollary A.19. Note that associated $\omega$-divisible module has strict $A$ action since it comes from division points of a divisible formal group. \\
One would like to show that $V$ is absolutely simple $\mathfrak{h}$ module. Note that in this case, $n_0 = h_r - 1$ and $n_1 = 1$ which are clearly relatively prime if $h_r \geq 2$. \\
Proof of this property is similar to  proof of \cite[Section-5, Proposition-8]{serre 2}. Here lie group is over $K(A)$. For this argument one needs a ramification theory result which goes over for $\omega$-torsion points and implies $G$ is infinite for all $h_r \geq 1$ (loc. cit. Lemma 3).\\
Now result follows from Lemma A.22. $\square$\\~\\
Proposition A.23 assumes lemma 3.7 in section-3 and proves proposition 3.8. Its content is same as proposition 3.9.\\
In light of developments in section-3 and this section we conclude that correct dimension implies $K(A)$ structure; $K(A)$-structure and Hodge-Tate decomposition implies correct dimension.

\vspace{1cm}
Soumyadip Sahu\\
soumyadip.sahu00@gmail.com

\end{document}